\documentclass[10pt]{birkjour}
\newtheorem{thm}{Theorem}[section]
\newtheorem{cor}[thm]{Corollary}
\newtheorem{lem}[thm]{Lemma}
\newtheorem{prop}[thm]{Proposition}
\theoremstyle{definition}

\theoremstyle{remark}

\numberwithin{equation}{section}
\begin{document}
	
	%
	%
	%
	%
	%

	%
	\author[Omar Ajebbar and Elhoucien Elqorachi]{Omar Ajebbar$^{1}$ and Elhoucien Elqorachi$^{2}$}
	\address{%
	$^1$
	Sultan Moulay Slimane University, Polydisciplinary Faculty,
	Department of Mathematics and Computer Science,
	Beni Mellal,
	Morocco}
	\address{%
		$^2$
	 Ibn Zohr University, Faculty of Sciences,
		Department of Mathematics,\\
		Agadir,
		Morocco}

\email{omar-ajb@hotmail.com, elqorachi@hotmail.com }

\title[Cosine-sine functional equation on semigroups]{Cosine-sine functional equation on semigroups}
\begin{abstract} Let $S$ be a semigroup. We determine the complex-valued solutions $f,g,h$ of the functional equation
\begin{equation*}f(xy)=f(x)g(y)+g(x)f(y)+h(x)h(y), x,y\in S,\end{equation*}
in terms of multiplicative functions, solutions of the special case
$$\varphi(xy)=\varphi(x)\chi(y)+\chi(x)\varphi(y), x,y\in S$$
of the sine addition law, where $\chi:S\to\mathbb{C}$ is a multiplicative function, and also in terms of solutions of the particular case
$$\psi(xy)=\psi(x)\chi(y)+\chi(x)\psi(y)+\varphi(x)\varphi(y), x,y\in S$$
of the cosine-sine functional equation where $\chi:S\to\mathbb{C}$ is a multiplicative function and $\varphi:S\to\mathbb{C}$ such that the pair $(\varphi,\chi)$ satisfies
the sine addition law.
\end{abstract}
\maketitle
\section{Introduction}Let $S$ be a semigroup. The functional equation
\begin{equation}\label{eq1}
f(xy)=f(x)g(y)+g(x)f(y)+h(x)h(y),\quad x,y\in S
\end{equation}
where $f,g,h:S\to\mathbb{C}$ are three unknown functions, is called the cosine-sine functional equation. It contained a number of fundamental functional equations which have been treated by several authors in the literature. For $h=0$ and $g=1$ we obtain the additive Cauchy equation. The special case of (\ref{eq1}) in which $h=0$ and $g=\dfrac{1}{2}f$ is the multiplicative Cauchy equation. For $h=0$ we obtain the sine addition law
\begin{equation*}
f(xy)=f(x)g(y)+g(x)f(y),\quad x,y\in S,
\end{equation*} which has been treated on groups, semigroups and algebras. See for example \cite[chapter 4]{H.St} and
\cite{EB-St}.
\par Chung, Kannappan and Ng \cite{CKN} solved
the functional (\ref{eq1}) on groups. They expressed the solutions as linear combination of additive, multiplicative functions and their product whose coefficients enter into matrix relations.\\
We refer also to \cite{Acz}, for further contextual and historical discussions.
 \par In \cite{AjElq} the functional equation has been extended to semigroups generated by their squares  by given more explicit solutions with simple constraints.
\par In \cite{H.St1} Stetk{\ae}r solved the functional equation \begin{equation}\label{new}
f(xy)=f(x)g(y)+g(x)f(y)-g(x)g(y), \;\;x,y\in S,
\end{equation}
where $S$ a semigroup.
\par Ebanks solved the functional equation (\ref{eq1}) on a t-compatible topological semigroup \cite [Theorem 5.1]{EB1}.
\par We are motivated to study the functional (\ref{eq1}) on a semigroup because it occurs as corollary in \cite [Corollary 6.1]{EB} in which $S$ is a topological commutative monoid.
\par Our work is in the same spirit as those of \cite{H.St1} and \cite{EB}. The contributions of the present parer to the knowledge about solutions of the functional equation (\ref{eq1}) are the following:
\begin{enumerate}
\item The setting has $S$ to be a general semigroup.
\item We relate the solutions of (\ref{eq1}) to those of the special case $\varphi(xy)=\varphi(x)\chi(y)+\chi(x)\varphi(y), x,y\in S$ of the sine addition law, where $\chi:S\to \mathbb{C}$ is a multiplicative function, and also to those of the functional equation $\psi(xy)=\psi(x)\chi(y)+\chi(x)\psi(y)+\varphi(x)\varphi(y), x,y\in S$, where $\psi:S\to \mathbb{C}$, $\chi:S\to \mathbb{C}$ is a multiplicative function such that the pair $(\varphi,\chi)$ satisfies the sine addition law.
\end{enumerate}
The moment we are submitting our work, an Ebanks paper \cite{new} has just appeared dealing with the same problem. Our paper does present the formulas for the solutions differently than \cite{new}, and it appears that the basic ideas in the proofs are very similar. Moreover, in our manuscript we have the presence of the explicit solutions of (\ref{new}) as an application of our main result.

\par The organization of the paper is as follows. In Sect. 2, we give notations and terminology. In Sect. 3, we give preliminary results that we will need in the paper. In Sect. 4 we prove our main results. In Sect. 5 we give an application of main results of the paper.
\section{Notations and terminology}
Throughout this paper $S$ denotes a semigroup. That is a set with an
associative composition.\\
A multiplicative function on $S$ is a homomorphism
$\chi:S\to(\mathbb{C},\cdot)$.\\
A function $f:S\to \mathbb{C}$ is said to be non-zero if there exists an element $x\in S$ such that $f(x)\neq0$. We denote $f\neq0$ a non-zero function.
\section{Preliminary results}
\begin{lem}\label{lem00}
Let $S$ be  a semigroup, $n\in\mathbb{N}\setminus\{0\}$, and $\chi_1,\chi_2,...\chi_n:S\to\mathbb{C}$ be different non-zero multiplicative functions. Then  \\
$\{\chi_1,\chi_2,...\chi_n\}$ is linearly independent.
\end{lem}
\begin{proof} See \cite [Theorem 3.18(b)]{H.St}.
\end{proof}
\begin{lem}\label{lem0}Let $f,g,:S\to\mathbb{C}$ be a solution of the sine addition law $f(xy)=f(x)g(y)+g(x)f(y),\,x,y\in S$ such that $f\neq0$. Then There exist a constant $\lambda\in\mathbb{C}$ and two multiplicative functions $\chi_{1},\chi_{2}:S\to\mathbb{C}$ such that $$2\,\lambda f=\chi_{1}-\chi_{2}\,\,\text{and}\,\,g=\dfrac{\chi_{1}+\chi_{2}}{2}.$$
\end{lem}
\begin{proof} See \cite [Theorem 4.1(b)]{H.St}.
\end{proof}
\begin{lem}\label{lem1} Let $f,g,h:S\to\mathbb{C}$ be a solution of (\ref{eq1}) such that $f$ and $h$ are linearly independent. Then there exist constants $\gamma_{1},\gamma_{2},\delta_{3},\delta_{4},\eta,\lambda,\mu\in\mathbb{C}$ such that
\begin{equation}\label{eq2}\begin{split}&g(xy)=g(x)g(y)+\gamma_{1}f(x)f(y)+\gamma_{2}[f(x)h(y)+h(x)f(y)]\\&\,+\delta_{3}h(x)h(y),\end{split}\end{equation}
\begin{equation}\label{eq3}\begin{split}&h(xy)=g(x)h(y)+h(x)g(y)+\gamma_{2}f(x)f(y)+\delta_{3}[f(x)h(y)+h(x)f(y)]\\&\,\,\,\,+\delta_{4}\,h(x)h(y),\,x,y\in S,\end{split}\end{equation}
\begin{equation}\label{eq4}(\lambda f+g)(xy)=(\lambda f+g)(x)(\lambda f+g)(y)+(\mu f+\eta h)(x)(\mu f+\eta h)(y)\end{equation}
for all $x,y\in S,$\\
\begin{equation}\label{eq5}\gamma_{1}+\gamma_{2}\delta_{4}-\delta_{3}^{2}=0\end{equation}
and
\begin{equation}\label{eq6}\lambda^{2}+\mu^{2}-\gamma_{1}=0,\quad\mu\eta-\gamma_{2}=0,\quad\lambda-\eta^{2}+\delta_{3}=0.\end{equation}
\end{lem}
\begin{proof} By using similar computations to the ones in the proof of \cite[Section 3, Theorem]{CKN}.
\end{proof}
In  Proposition \ref{prop1} we show how the solutions are related to those of the following special case of the sine addition law:
\begin{equation*}\varphi(xy)=\varphi(x)\chi(y)+\chi(x)\varphi(y), x,y\in S,\end{equation*}
where $\chi:S\to\mathbb{C}$ is a multiplicative function.
\begin{prop}\label{prop1} Let $f,g,h:S\to\mathbb{C}$ be solutions of (\ref{eq1}) such that $f$ and $h$ are linearly independent. Then the triplet $(f,g,h)$ satisfies one of the following possibilities:\\
(1) There exist a constant $\beta\in\mathbb{C}$ and two multiplicative functions $\mu,\chi:S\to\mathbb{C}$ such that
\begin{equation}\label{eq06}g=\chi,\,h-\beta f=\varphi\,\,\text{and}\,\,\chi+\beta h=\mu,\end{equation}
where $\varphi:S\to\mathbb{C}$ satisfies the special case of the sine addition law  $\varphi(xy)=\varphi(x)\chi(y)+\chi(x)\varphi(y), x,y\in S$.\\
(2) There exist constants $\alpha\in\mathbb{C}\setminus\{0\}, \beta\in\mathbb{C}$ and a multiplicative function $\chi:S\to\mathbb{C}$ such that
\begin{equation}\label{eq006}-\alpha\,f+g=\chi\end{equation}
and the pair $(h,\chi+2\,\alpha f+\dfrac{1}{2}\,\beta h)$ satisfies the sine addition law, i.e.,
\begin{equation}\label{eq0006}h(xy)=h(x)(\chi+2\,\alpha f+\dfrac{1}{2}\,\beta h)(y)+(\chi+2\,\alpha f+\dfrac{1}{2}\,\beta h)(x)h(y),\,x,y\in S.\end{equation}
(3) $$f=F,\,\,\,\,g=-\dfrac{1}{2}\delta^{2}\,F+G+\delta\,H,\,\,\,\,h=-\delta
\,F+H,$$
where $\delta\in\mathbb{C}$ is a constant and $F,G,H:S\to\mathbb{C}$ are functions such that the triplet $(F,G,H)$ satisfies (1)-(2).
\end{prop}
\begin{proof}Let $f,g,h:S\to\mathbb{C}$ be a solution of (\ref{eq1}) such that $f$ and $h$ are linearly independent. Then, the functional equations (\ref{eq2}), (\ref{eq3}) and (\ref{eq4}), and the identities (\ref{eq5}) and (\ref{eq6}) are satisfied. We split the discussion into the cases $\gamma_{2}=0$ and $\gamma_{2}\neq0$ \\
\underline{Case A}: $\gamma_{2}=0$. Then, we get from (\ref{eq5}) that $\gamma_{1}=\delta_{3}^{2}$. Notice that we can choose $\mu=\eta=0$ in (\ref{eq6}). So that $\lambda=-\delta_{3}$. Hence, the functional equation (\ref{eq4}) reduces to
\begin{equation}\label{eq7}(-\delta_{3}f+g)(xy)=(-\delta_{3}f+g)(x)(-\delta_{3}f+g)(y)\end{equation}
for all $x,y\in S,$ which implies that there exists a multiplicative function $\chi:S\to\mathbb{C}$ such that
\begin{equation}\label{eq8}-\delta_{3}f+g=\chi.\end{equation}
When we substitute this in  (\ref{eq1}) and  (\ref{eq3}) we get respectively the following by a small computation
\begin{equation}\label{eq9}f(xy)=f(x)\chi(y)+\chi(x)f(y)+2\delta_{3}f(x)f(y)+h(x)h(y)\end{equation}
and
\begin{equation}\label{eq10}h(xy)=h(x)(\chi+2\delta_{3}f+\dfrac{1}{2}\delta_{4}h)(y)+(\chi+2\delta_{3}f+\dfrac{1}{2}\delta_{4}h)(x)h(y)\end{equation}
for all $x,y\in S.$
We have the following subcases.\\
\underline{Subcase A.1}: $\delta_{3}=0$. Then, we get from (\ref{eq8}) that $g=\chi$. So, (\ref{eq9}) and (\ref{eq10}) reduce respectively to
\begin{equation}\label{eq11}f(xy)=f(x)\chi(y)+\chi(x)f(y)+h(x)h(y)\end{equation}
and
\begin{equation}\label{eq12}h(xy)=h(x)\chi(y)+\chi(x)h(y)+\delta_{4}h(x)h(y)\end{equation}
for all $x,y\in S.$ So that $$h(xy)-\delta_{4}f(xy)=h(x)\chi(y)+\chi(x)h(y)-\delta_{4}f(x)\chi(y)-\delta_{4}\chi(x)f(y),$$ for all $x,y\in S$, which implies that \begin{equation}\label{eq13}(h-\delta_{4}f)(xy)=(h-\delta_{4}f)(x)\chi(y)+\chi(x)(h-\delta_{4}f)(y),\end{equation} for all $x,y\in S.$
Now, using  (\ref{eq12}) we obtain the functional equation
\begin{equation*}(\chi+\delta_{4}h)(xy)=(\chi+\delta_{4}h)(x)(\chi+\delta_{4}h)(y),\end{equation*}
for all $x,y\in S$. Hence, the function $\chi+\delta_{4}\,h$ is multiplicative. So, by writing $\beta$ instead of $\delta_{4}$ and putting $\varphi:=h-\delta_{4}f$ in (\ref{eq13}) and $\mu:=\chi+\delta_{4}h$ we get part (1) in Proposition \ref{prop1}.\\
\underline{Subcase A.1}: $\delta_{3}\neq0$. Then by putting $\alpha:=\delta_{3}$ and $\beta:=\delta_{4}$ and applying Lemma \ref{lem0} to  (\ref{eq10}) we get part (2).\\
\underline{Case B}: $\gamma_{2}\neq0$. By proceeding as in \cite[Section 3, Case 2.]{CKN} we obtain part (3).
\end{proof}
\begin{prop}\label{prop01} Let $\chi:S\to\mathbb{C}$ be a multiplicative function and $\varphi:S\to\mathbb{C}$  a function such that $(\varphi,\chi)$ satisfies the sine addition law $$\varphi(xy)=\varphi(x)\chi(y)+\chi(x)\varphi(y),\,x,y\in S$$
If $\varphi=a\chi_{1}+b\chi_{2}$, where $a,b\in \mathbb{C}$ and $\chi_{1},\chi_{2}:S\to\mathbb{C}$ non-zero different multiplicative functions, then $\varphi=0$.
\end{prop}
\begin{prop}\label{prop2} Let $\chi:S\to\mathbb{C}$ be a multiplicative function and $\psi,\varphi:S\to\mathbb{C}$ tow non-zero functions such that $(\varphi,\chi)$ satisfies the sine addition law $$\varphi(xy)=\varphi(x)\chi(y)+\chi(x)\varphi(y),\,x,y\in S$$ and $(\psi,\chi,\varphi)$ satisfies the cosine-sine functional equation
$$\psi(xy)=\psi(x)\chi(y)+\chi(x)\psi(y)+\varphi(x)\varphi(y),\,x,y\in S.$$
If $\alpha,\beta,\gamma\in\mathbb{C}$ such that $\alpha \psi+\beta \chi+\gamma \varphi=0$ then $\alpha=\gamma=0$ and $\beta \chi=0$. In other words if $\chi\neq0$ then $\psi,\,\chi$ and $\varphi$ are linearly independent.
\end{prop}
\section{Main results}
In this section we describe  the solutions of (\ref{eq1}) in a semigroup $S$ in terms of multiplicative functions on $S$, solutions
 $\varphi(xy)=\varphi(x)\chi(y)+\chi(x)\varphi(y), x,y\in S$ of the sine addition law, where $\chi:S\to \mathbb{C}$ is a multiplicative function, and also solutions of the particular functional equation $\psi(xy)=\psi(x)\chi(y)+\chi(x)\psi(y)+\varphi(x)\varphi(y), x,y\in S$ such that the pair $(\varphi,\chi)$ satisfies the sine addition law.\\
 The method of proof parallels \cite{CKN} and \cite{EB1}.
\subsection{Solutions of (\ref{eq1}) when $f$ and $h$ are linearly dependent}
\begin{thm}\label{thm1} The solutions $f,g,h:S\to\mathbb{C}$ of the functional equation (\ref{eq1}) such that $f$ and $h$ are linearly dependent are one of the following forms\\
(1) $f=0$, $g$ arbitrary and $h=0$.\\
(2) $f=\lambda\,(\chi_{1}-\chi_{2})$,
$g=\frac{\chi_{1}+\chi_{2}}{2}-\rho^{2}\,\frac{\chi_{1}-\chi_{2}}{2\lambda}$ and
$h=\rho\,(\chi_{1}-\chi_{2})$ where
$\lambda\in\mathbb{C}\setminus\{0\}$ and $\rho\in\mathbb{C}$ are two
constants, $\chi_{1},\,\chi_{2}:S\to\mathbb{C}$ are two
multiplicative functions such that $\chi_{1}\neq\chi_{2}$.\\
(3) $f=\varphi,\,\,g=\chi-\frac{c^{2}}{2}\,\varphi\,\,\text{and}\,\,h=c\,\varphi$, where $c\in\mathbb{C}$ is a constant, $\chi:S\to\mathbb{C}$ a
multiplicative function and $\varphi:S\to\mathbb{C}$ a non-zero solution of the special case of the sine addition law $\varphi(xy)=\varphi(x)\chi(y)+\chi(x)\varphi(y), x,y\in S$.
\end{thm}
\begin{proof} Let $f,g,h:S\to\mathbb{C}$ be a solution of
the functional equation (\ref{eq1}) such that $f$ and $h$ are linearly
dependent. If $f=0$, then (\ref{eq1}) implies
$h(x)h(y)=0$ for all $x,y\in S$, so $h=0$ and $g$ is arbitrary.
\\In what remains of the proof we assume that $f\neq0$. Since $f$ and $h$ are linearly
dependent there exists a constant $c\in\mathbb{C}$ such that $h=cf$. So, (\ref{eq1}) becomes
\begin{equation}\label{eq012}
f(xy)=f(x)k(y)+k(x)f(y), x,y\in S,
\end{equation}
where $k:=g+\frac{c^{2}}{2}\,f$. \\Hence, according to Lemma \ref{lem0} we get that
 $2\alpha f=\chi_{1}-\chi_{2}$ and $k=\frac{\chi_{1}+\chi_{2}}{2}$ where
$\alpha\in\mathbb{C}$ is a constant and
$\chi_{1},\,\chi_{2}:S\to\mathbb{C}$ are two multiplicative
functions.\\We have two cases to consider.\\
\underline{Case A}: $\alpha\neq0$. So, by putting $\lambda:=\frac{1}{2\alpha}$ and $\rho:=\lambda c$ we get part (2) as in Case 1 of the proof of \cite [Theorem 4.2]{AjElq}.\\
\underline{Case B}: $\alpha=0$. Then $\chi_{1}=\chi_{2}$. So $k=\chi$ and from (\ref{eq012}) we get that $f=\varphi$, where $\chi:=\chi_{1}$ and $\varphi$ is a non-zero solution of the special case of the sine addition law $\varphi(xy)=\varphi(x)\chi(y)+\chi(x)\varphi(y), x,y\in S$. Hence $g=\chi-\frac{c^{2}}{2}\,\varphi$ and $h=c\,\varphi$.
\par Conversely, if $f$, $g$ and $h$ are of the forms (1)-(3) in Theorem
\ref{thm1} we check by elementary computations that $f$, $g$ and $h$
satisfy the functional equation (\ref{eq1}), and that $f$ and $h$ are
linearly dependent. This completes the proof of Theorem \ref{thm1}.
\end{proof}
\subsection{Solutions of (\ref{eq1}) when $f$ and $h$ are linearly independent}
\begin{thm}\label{thm2} The solutions $f,g,h:S\to\mathbb{C}$ of the functional equation (\ref{eq1}) such that $f$ and $h$ are linearly independent are of form
\begin{center}
$\begin{pmatrix} f \\
g \\ h \end{pmatrix}=\begin{pmatrix} 1& 0 & 0 \\
-\frac{\delta^{2}}{2} & 1 & \delta \\ -\delta & 0 & 1 \\ \end{pmatrix}\begin{pmatrix} F\\
G \\ H \end{pmatrix}$,\end{center}
where $\delta\in\mathbb{C}$ is a constant and the functions
$F,G,H:S\to\mathbb{C}$ are of the following families:\\
(1)
\begin{center}
$\begin{pmatrix} F \\
G \\ H \end{pmatrix}=\begin{pmatrix} 1 & 0 & 0 \\
0 & 1 & 0 \\ 0 & 0 & 1 \end{pmatrix}\begin{pmatrix} \psi\\
\chi \\ \varphi \end{pmatrix}$,\end{center}
where  $\chi:S\to\mathbb{C}$ is multiplicative function, $\psi,\varphi:S\to\mathbb{C}$ are non-zero functions such that $(\varphi,\chi)$ satisfies the sine addition law  $\varphi(xy)=\varphi(x)\chi(y)+\chi(x)\varphi(y), x,y\in S$ and the triplet $(\psi,\chi,\varphi)$ satisfies the cosine-sine functional equation $\psi(xy)=\psi(x)\chi(y)+\chi(x)\psi(y)+\varphi(x)\varphi(y), x,y\in S$.\\
(2)
\begin{center}
$\begin{pmatrix} F \\
G \\ H \end{pmatrix}=\begin{pmatrix} c^{2} & -c^{2} & -c \\
0 & 1 & 0 \\ c & -c & 0 \\ \end{pmatrix}\begin{pmatrix} \mu\\
\chi \\ \varphi \end{pmatrix}$,\end{center}
where $c\in\mathbb{C}\setminus\{0\}$ is a constant,
$\chi,\,\mu:S\rightarrow\mathbb{C}$ are two multiplicative functions
with $\chi\neq\mu$, and $\varphi:S\to\mathbb{C}$ is a non-zero function such that $(\varphi,\chi)$ satisfies the sine addition law  $\varphi(xy)=\varphi(x)\chi(y)+\chi(x)\varphi(y), x,y\in S$.\\
(3)
\begin{center}
$\begin{pmatrix} F \\
G \\ H \end{pmatrix}=\begin{pmatrix} -c_{1} & c_{1} & -c_{1}c_{2} \\
\frac{1}{2} & \frac{1}{2} & -\frac{1}{2}c_{2} \\ 0 & 0 & 1\\ \end{pmatrix}\begin{pmatrix} \mu\\
\chi \\ \varphi  \end{pmatrix}$, \end{center}
where $c_{1},c_{2}\in\mathbb{C}\setminus\{0\}$ are two constants
satisfying $1+c_{1}\,c_{2}^{2}=0$, $\chi,\,\mu:S\rightarrow\mathbb{C}$ are two multiplicative functions
such that $\chi\neq\mu$, and $\varphi:S\to\mathbb{C}$ is a non-zero function such that $(\varphi,\chi)$ satisfies the sine addition law  $\varphi(xy)=\varphi(x)\chi(y)+\chi(x)\varphi(y), x,y\in S$.\\
(4)
\begin{center}
$\begin{pmatrix} F \\
G \\ H \end{pmatrix}=\begin{pmatrix} c\,\rho & c\,(2-\rho) & -2\,c \\
\frac{1}{4}\,\rho & \frac{1}{4}\,(2-\rho) & \frac{1}{2} \\ \frac{1}{2\,\lambda} & -\frac{1}{2\,\lambda} & 0 \\ \end{pmatrix}\begin{pmatrix} \chi_{1}\\
\chi_{2} \\ \chi_{3} \end{pmatrix}$,\end{center} where $\lambda,
\rho,c\in\mathbb{C}\setminus\{0\}$ are three constants with
$2\,c\,\lambda^{2}\,\rho\,(2-\rho)=1$;
$\chi_{1},\,\chi_{2},\,\chi_{3}:S\rightarrow\mathbb{C}$ are three
multiplicative functions such that $\chi_{1}\neq\chi_{2}$,
$\chi_{1}\neq\chi_{3}$, $\chi_{2}\neq\chi_{3}$.\\
\end{thm}
\begin{proof} According to Proposition \ref{prop1} we have one of the following cases:\\
\underline{Case A}: $(f,g,h)$ is of the form in Proposition \ref{prop1}(1).
We discuss according to whether $\beta=0$ or $\beta\neq0$.\\
\underline{Subcase A.1}: $\beta=0$. Then we get from (\ref{eq06}) that $g=\chi$ and $h=\varphi$.  So, $f$ satisfies the special case the cosine-sine functional equation $$f(xy)=f(x)\chi(y)+\chi(x)f(y)+\varphi(x)\varphi(y), x,y\in S.$$ Notice that $\varphi\neq0$ and $f\neq0$ because $f$ and $h$ are linearly independent. The result occurs in part (1) $\delta=0$.\\
\underline{Subcase A.2}: $\beta\neq0$. By putting $c:=\dfrac{1}{\beta}$ we derive from  (\ref{eq06}) that $g=\chi$, $f=c^{2}\mu-c^{2}\chi-c\,\varphi$ and $h=c\,\mu-c\,\chi$. The linear independence of $f$ and $h$ imposes the conditions $c\neq0$ and $\chi\neq\mu$. So obtain part (2) $\delta=0$.\\
\underline{Case B}: $(f,g,h)$ is of the form in Proposition \ref{prop1}(2). Then, according to Lemma \ref{lem0}, there exist a constant $\lambda\in\mathbb{C}$ and two multiplicative functions $\chi_{1},\chi_{2}:S\to\mathbb{C}$ such that
\begin{equation}\label{eq14} 2\,\lambda h=\chi_{1}-\chi_{2}\end{equation}
and
\begin{equation}\label{eq15}\chi+2\,\alpha f+\dfrac{1}{2}\,\beta h=\dfrac{\chi_{1}+\chi_{2}}{2}.\end{equation}
As in Case A we discuss according to whether $\lambda=0$ or $\lambda\neq0$.\\
\underline{Subcase B.1}: $\lambda=0$. Then, we deduce from (\ref{eq14}) that $\chi_{1}=\chi_{2}$. So, in view of (\ref{eq15}) and (\ref{eq0006}) and putting $\mu:=\chi_{1}=\chi_{2}$, we derive that there exists a solution $\varphi:S\to\mathbb{R}$ of the special case of the sine addition law $\varphi(xy)=\varphi(x)\mu(y)+\mu(x)\varphi(y),x,y\in S$ such that
\begin{equation}\label{eq16}h=\varphi.\end{equation} As $f$ and $h$ are linearly independent we get that $\varphi\neq0$ and $\chi\neq\mu$. Now, with $c_{1}:=\dfrac{1}{2\alpha }$ and $c_{2}:=\dfrac{1}{2}\beta$ we get from (\ref{eq15}) and (\ref{eq006}) that
\begin{equation}\label{eq17}f=-c_{1}\chi+c_{1}\mu-c_{1}c_{2}\varphi\end{equation} and
\begin{equation}\label{eq18}g=\dfrac{1}{2}\chi+\dfrac{1}{2}\mu-\dfrac{1}{2}c_{2}\varphi.\end{equation} Moreover, using (\ref{eq16}), (\ref{eq17}) and (\ref{eq18}), a small computation shows that
\begin{equation*}\begin{split}&f(x)g(y)+g(x)f(y)+h(x)h(y)\\
&=[-c_{1}\chi(x)+c_{1}\mu(x)-c_{1}c_{2}\varphi(x)][\frac{1}{2}\chi(y)+\frac{1}{2}\mu(y)-\frac{1}{2}c_{2}\varphi(y)]\\
&+[-c_{1}\chi(y)+c_{1}\mu(y)-c_{1}c_{2}\varphi(y)][\frac{1}{2}\chi(x)+\frac{1}{2}\mu(x)-\frac{1}{2}c_{2}\varphi(x)]+\varphi(x)\varphi(y)\\
&=f(xy)+(1+c_{1}c_{2}^{2})\varphi(x)\varphi(y),
\end{split}\end{equation*}
for all $x,y\in S$. Since $(f,g,h)$ is a solution of (\ref{eq1})  and $\varphi\neq0$ we derive from the identity above that $1+c_{1}c_{2}^{2}=0$. Hence, by interchanging  $\mu$ and $\chi$, we obtain a solution of the family (3) for $\delta=0$.\\
\underline{Subcase B.2}: $\lambda\neq0$. Then, from (\ref{eq14}) we get that
\begin{equation}\label{eq19}h=\frac{\chi_{1}-\chi_{2}}{2\lambda}\end{equation}
Now, by putting $c:=\dfrac{1}{4\alpha}$ and $\rho:=1-\dfrac{\beta}{2\lambda}$ we drive from (\ref{eq15}) and (\ref{eq006}) respectively
\begin{equation}\label{eq20}f=c\rho\chi_{1}+c(2-\rho)\chi_{2}-2c\chi\end{equation}
and
\begin{equation}\label{eq21}g=\frac{1}{4}\rho\,\chi_{1}+\frac{1}{4}(2-\rho)\chi_{2}+\frac{1}{2}\chi.\end{equation}
On the other hand, by using (\ref{eq19}), (\ref{eq20}) and (\ref{eq21}) we obtain
\begin{equation*}\begin{split}
&f(x)g(y)+g(x)f(y)+h(x)h(y)\\
&=(\frac{c\rho^{2}}{2}+\frac{1}{4\rho^{2}})\chi_{1}(xy)+(\frac{c\rho(2-\rho)}{2}-\frac{1}{4\lambda^{2}})\chi_{1}(x)\chi_{2}(y)\\
&\quad+(\frac{c\rho(2-\rho)}{2}-\frac{1}{4\lambda^{2}})\chi_{1}(y)\chi_{2}(x)+(\frac{c(2-\rho)^{2}}{2}+\frac{1}{4\lambda^{2}})\chi_{2}(xy)-2c\chi(xy)\\
&=[(\frac{c\rho^{2}}{2}+\frac{1}{4\lambda^{2}})\chi_{1}(y)+(\frac{c\rho(2-\rho)}{2}-\frac{1}{4\lambda^{2}})\chi_{2}(y)]\chi_{1}(x)\\
&\quad+[(\frac{c\rho(2-\rho)}{2}-\frac{1}{4\lambda^{2}})\chi_{1}(y)+(\frac{c(2-\rho)^{2}}{2}+\frac{1}{4\lambda^{2}})\chi_{2}(y)]\chi_{2}(x)-2c\chi(y)\chi(x)
\end{split}\end{equation*}
and
$$f(xy)=c\rho\chi_{1}(y)\chi_{1}(x)+c(2-\rho)\chi_{2}(y)\chi_{2}(x)-2c\chi(y)\chi(x),$$
for all $x,y\in S$.\\
As $\chi_{1}$ et $\chi_{2}$ are different multiplicative functions we derive from the two identities above and (\ref{eq1}),  by applying Lemma \ref{lem00}, that we have at least one of the following identities ( according to $\chi_{1}=0$ and $\chi_{2}\neq0$, or $\chi_{2}=0$ and $\chi_{1}\neq0$, or $\chi_{1}\neq0$ and $\chi_{2}\neq0$):
\[
\left\{
\begin{array}{r c l}
\frac{c\rho^{2}}{2}+\frac{1}{4\lambda^{2}}=\rho c\\
\frac{c\rho(2-\rho)}{2}-\frac{1}{4\lambda^{2}}=0\\
\frac{c(2-\rho)^{2}}{2}+\frac{1}{4\rho^{2}}=c(2-\rho),
\end{array}
\right.
\]
Hence, $2c\lambda^{2}\rho(2-\rho)=1$. Moreover, since $f$ and $h$ are linearly independent we derive from (\ref{eq19}) and (\ref{eq20}) that $\chi_{1}\neq\chi$ and $\chi_{2}\neq\chi$. So, By putting $\chi_{3}:=\chi$ we get the solution of part (4) $\delta=0$.\\
\underline{Case C}: $(f,g,h)$ is of the form in Proposition \ref{prop1}(3). So, we get one of the parts (1)-(4) for an arbitrary constant $\delta\in\mathbb{C}$.
\par Conversely, if $(f,g,h)$ is  of the form in Theorem \ref{thm2} such that $(F,G,H)$ is of the forms (1)-(4) in Theorem \ref{thm2} we check by elementary computations that $(f,g,h)$
satisfies the functional equation (\ref{eq1}), and that $f$ and $h$ are
linearly independent. This completes the proof of Theorem \ref{thm2}.
\end{proof}
\section{Application}
In \cite {H.St1} a special case of (\ref{eq1}), in which $h=i g$, i.e.,
\begin{equation}\label{eq22}f(xy)=f(x)g(y)+g(x)f(y)-g(x)g(y), x,y\in S\end{equation}
where $S$ is a semigroup, was solved.\\
Now, we find the result in \cite [Theorem 3.3]{H.St1} in Corollary \ref{cor1} as a consequence of Theorem \ref{thm1} and Theorem \ref{thm2}.
\begin{cor}\label{cor1}The solutions $f,g:S\to\mathbb{C}$ of the functional equation (\ref{eq22}) are one of the following forms\\
(1) $f$ is a function such that $f(xy)=0$ for all $x,y\in S$ and $g=0$.\\
(2) $f$ is a non-zero function such that $f(xy)=0$ for all $x,y\in S$ and $g=2 f$.\\
(3) $$f=\frac{\alpha^{2}}{2\alpha-1}\chi\,\,\text{and}\,\,g=\alpha \chi,$$
where $\chi:S\to \mathbb{C}$ is a non-zero multiplicative function and $\alpha\in\mathbb{C}\setminus\{0,\frac{1}{2}\}$ is a constant.\\
(4) $$f=\frac{\chi_{1}+\chi_{2}}{2}+\frac{\beta^{2}+1}{2 \beta}\frac{\chi_{1}-\chi_{2}}{2}\,\,\text{and}\,\,g=\frac{\chi_{1}+\chi_{2}}{2}+\beta\frac{\chi_{1}-\chi_{2}}{2},$$
where $\chi_{1},\chi_{2}:S\to \mathbb{C}$ are two different multiplicative functions and $\beta\in\mathbb{C}\setminus\{0\}$ is a constant.\\
(5) $$f=\frac{1}{2}\varphi+\chi\,\,\text{and}\,\,g=\varphi+\chi,$$
where $\chi:S\to \mathbb{C}$ is a non-zero multiplicative function and $\varphi:S\to \mathbb{C}$ a non-zero solution of sine addition law $\varphi(xy)=\varphi(x)\chi(y)+\chi(x)\varphi(y), x,y\in S$.\\
(6) $$f=\varphi+\chi\,\,\text{and}\,\,g=\chi,$$
where $\chi:S\to \mathbb{C}$ is a non-zero multiplicative function and $\varphi:S\to \mathbb{C}$ a non-zero solution of sine addition law $\varphi(xy)=\varphi(x)\chi(y)+\chi(x)\varphi(y), x,y\in S$.
\end{cor}
\begin{proof} The solutions $f,g:S\to\mathbb{C}$ of (\ref{eq22}) are such that the triplet $(f,g,ig)$ is a solution of (\ref{eq1}). So, we discuss according to whether $f$ and $g$ are linearly dependent or $f$ and $g$ are linearly independent.\\
\underline{Case A}: $f$ and $g$ are linearly dependent. According to Theorem \ref{thm1} we have the following possibilities:
\begin{enumerate}
\item $f=0$ and $g=0$. This is a special case of part (1).
\item The triplet $(f,g,i g)$ is of the form in Theorem \ref{thm1}(2). Then
\begin{equation}\label{eq24}(\rho-\frac{i}{2}+\frac{i \rho^{2}}{2 \lambda})\chi_{1}+(\rho+\frac{i}{2}+\frac{i \rho^{2}}{2 \lambda})\chi_{2}=0.\end{equation} As the coefficients $\rho-\frac{i}{2}+\frac{i \rho^{2}}{2 \lambda}$ and $\rho+\frac{i}{2}+\frac{i \rho^{2}}{2 \lambda}$ can not be zero at the same time we get from (\ref{eq24}), according to Lemma \ref{lem00}, that $\chi_{1}=0$ or $\chi_{2}=0$. Without loss of generality we assume $\chi_{1}=0$, then we get from Theorem \ref{thm1}(2) that
\begin{equation}\label{eq25}f=-\lambda \chi_{2}\quad\text{and}\quad g=\frac{\lambda+\rho^{2}}{2 \lambda}=i \rho \chi_{2}.\end{equation} Moreover, since $0=\chi_{1}\neq\chi_{2}$ we derive from (\ref{eq24}) that  $\rho+\frac{i}{2}+\frac{i \rho^{2}}{2 \lambda}=0$ which implies that $\rho\neq0$ because $i\neq0$ and  $\alpha^{2}=-\lambda(2 \alpha-1)$ with $\alpha:=i\rho$. So, (\ref{eq25}) can be rewritten as follows $$f=\frac{\alpha^{2}}{2\alpha-1}\chi\,\,\text{and}\,\,g=\alpha \chi,$$ with $\chi:=\chi_{2}$ and $\alpha\in\mathbb{C}\setminus\{0,\frac{1}{2}\}$. The solution is in (3).
\item The triplet $(f,g,i g)$ is of the form in Theorem \ref{thm1}(3).
Then, from the expressions of $g$ in Theorem \ref{thm1}(3), we get that
\begin{equation}\label{eq27}\frac{c(c-2i)}{2}\varphi=\chi.\end{equation} Since $(\varphi,\chi)$ satisfies the sine addition law  we get that $\frac{c(c-2i)}{2}\chi(xy)=0$ for all $x,y\in S$. So, for $\lambda^{2}:=\frac{c(c-2i)}{2}$ we obtain $\lambda \chi=0$. Hence $\lambda^{3} \varphi=0$. Then $\lambda=0$ because $\varphi\neq0$. It follows that $c=0$ or $c=2i$. So, in view of (\ref{eq27}) we get that $\chi=0$. So that $f(xy)=\varphi(xy)=\varphi(x)0+\varphi(y)0=0$ for all $x,y\in S$
\par If $c=0$ then $g=0$. The solution occurs in part (1).
\par If $c=2i$ then $g=-\frac{c^{2}}{2}\,\varphi=2 f$. The solution is in part (2).
\end{enumerate}
\underline{Case B}: $f$ and $g$ are linearly independent. According to Theorem \ref{thm2} we have
\begin{equation}\label{eq28} f=F, g=-\frac{1}{2}\delta^{2}F+G+\delta\,H\quad\text{and}\quad ig=-\delta\,F+H,\end{equation}
where $\delta\in\mathbb{C}$ is a constant and the functions
$F,G,H:S\to\mathbb{C}$ are of the forms (1)-(4) in Theorem \ref{thm2} with the same constraints.\\Then $-\delta\,F+H=-\frac{1}{2}i\delta^{2}F+iG+i\delta\,H$.
So that
\begin{equation}\label{eq29} \xi F=G+\eta H,\end{equation}
where
\begin{equation}\label{eq30} \xi:=i\delta+\frac{1}{2}\delta^{2}\end{equation}
and
\begin{equation}\label{eq31} \eta:=i+\delta.\end{equation}
We have the following possibilities for the triplet $(F,G,H)$.
\begin{enumerate}
\item The triplet $(F,G,H)$ is of the form in Theorem \ref{thm2}(1). Then, in view of (\ref{eq29}) we get that  $\xi\psi=\chi+\eta\varphi$. Applying Proposition \ref{prop2} we obtain $i\delta+\frac{1}{2}\delta^{2}=0$ and $i+\delta=0$, which is a contradiction. Hence this possibility does not arise.
\item The triplet $(F,G,H)$ is of the form in Theorem \ref{thm2}(2). Then, using (\ref{eq29}) we obtain
\begin{equation*}\xi c\varphi=(\xi c^{2}-\eta c)\mu-(\xi c^{2}-\eta c+1)\chi.\end{equation*}
Now, according to Proposition \ref{prop01} we derive from the identity above and (\ref{eq30}) that $\delta=0$ and $\eta=i$, or $\delta=-2i$ and $\eta=-i$ because $c\neq0$ and $\varphi\neq0$. So that
\begin{equation}\label{eq32}\pm ic \mu+(\mp ic+1)\chi=0.\end{equation}
Since $\mu,\chi:S\to\mathbb{C}$ are different multiplicative functions and $ic\neq0$ we deduce from (\ref{eq32}), according to Lemma \ref{lem00}, that $\mu=0$ and $\chi\neq0$.
\par Hence, $$\delta=0,\,\eta=i,\,\mu=0\quad\text{and}\quad c=-i$$
or $$\delta=-2i,\,\eta=-i,\,\mu=0\quad\text{and}\quad c=i.$$
Therefore $$f=F=\chi\pm i\varphi\quad\text{and}\quad g=\chi$$
The solution occurs in part (6) by writing $\varphi$ instead of $\pm i\varphi$.
\item The triplet $(F,G,H)$ is of the form in Theorem \ref{thm2}(3). Then, using (\ref{eq29})  and the formulas of $F,\,G$ and $H$ we derive
\begin{equation}\label{eq33}(\frac{1}{2}c_{2}-\xi\,c_{1}c_{2}-\eta)\varphi=(\frac{1}{2}+\xi\,c_{1})\mu+(\frac{1}{2}-\xi\,c_{1})\chi.\end{equation}
Since $\varphi\neq0$ we deduce from (\ref{eq33}), by applying Proposition \ref{prop01}, that
\begin{equation}\label{eq34}\frac{1}{2}c_{2}-\xi\,c_{1}c_{2}-\eta=0\end{equation}
and
\begin{equation}\label{eq35}(\frac{1}{2}+\xi\,c_{1})\mu+(\frac{1}{2}-\xi\,c_{1})\chi=0.\end{equation}
\par If $\chi=0$ then $\mu\neq0$ because $\chi\neq\mu$. So, (\ref{eq35}) and (\ref{eq33}) imply that $\xi\,c_{1}=-\frac{1}{2}$ and $\eta=c_{2}$. Hence, using (\ref{eq31}), $\delta=c_{2}-i$. Substituting this in (\ref{eq30}) and taking into account that  $\xi\,c_{1}=-\frac{1}{2}$, a small computation shows that $c_{1}=0$, which contradicts that $1+c_{1}\,c_{2}^{2}=0$.
\par Hence, since $\chi\neq\mu$ and the coefficients $\frac{1}{2}+\xi\,c_{1}$ and $\frac{1}{2}-\xi\,c_{1}$ can nHencet be zero at the same time we derive, by applying  Lemma \ref{lem00} to (\ref{eq35}), that $\mu=0$ and then $\chi\neq0$. So,in view of (\ref{eq34}) we get that $\xi\,c_{1}=\frac{1}{2}$. Combining this with (\ref{eq33}), (\ref{eq31})  and (\ref{eq30}) and taking the identity $1+c_{1}\,c_{2}^{2}=0$ into account, we deduce that $\eta=0$, $\delta=-i$, $\xi=\frac{1}{2}$, $c_{1}=1$ and $c_{2}=\pm i$. Hence, $F=\chi-c_{2} \varphi$, $G=\frac{1}{2}\chi-\frac{1}{2}c_{2} \varphi$ and $H=\varphi$. Now, by substituting this in (\ref{eq28}) we derive that $f=\chi-i\varphi$ and $g=\chi-(i+c_{2})\varphi$.
\par If $c_{2}=i$ the solution occurs in part (5).
\par If $c_{2}=-i$ the solution occurs in part (6).
\item The triplet $(F,G,H)$ is of the form in Theorem \ref{thm2}(4). According to (\ref{eq29}) we get that
\begin{equation}\label{eq36}(\xi c\rho-\frac{\rho}{4}-\frac{\eta}{2\lambda})\chi_{1}+(\xi c(2-\rho)-\frac{2-\rho}{4}+\frac{\eta}{2\lambda})\chi_{2}-(2\xi c+\frac{1}{2})\chi_{3}=0.\end{equation}
As the coefficients in (\ref{eq36}) can be zero at the same time, and the multiplicative functions $\chi_{1},\chi_{2}$ and $\chi_{3}$ are different we deduce, according to Lemma \ref{lem00} that just one of the multiplicative functions $\chi_{1},\chi_{2}$ and $\chi_{3}$ is zero. So have the following possibilities:
\begin{enumerate}
\item $\chi_{1}=0$. Then  $\chi_{2}\neq0$ and $\chi_{3}\neq$. As $\chi_{2}\neq\chi_{3}$ we derive from (\ref{eq36}), using Lemma \ref{lem00}, that
\begin{equation*}2\xi c+\frac{1}{2}=0\end{equation*}
and
\begin{equation*}2\xi c-\xi c\rho+\frac{\rho}{4}+\frac{\eta}{2\lambda}-\frac{1}{2}=0.\end{equation*}
Hence,
\begin{equation}\label{eq37}4\xi c=-1\end{equation}
and
\begin{equation}\label{eq38}\eta=\lambda(2-\rho).\end{equation}
Then, by using (\ref{eq31}) and (\ref{eq38}) and taking into account that
$2c\lambda^{2}\rho(2-\rho)=1$, we get that
\begin{equation}\label{eq39}\delta=\lambda(2-\rho)-i=-\frac{1+2i\lambda}{2\lambda},\end{equation} and then , after a small computation, that
\begin{equation}\label{eq40}c\rho=-1.\end{equation}
So that
\begin{equation}\label{eq41}2\lambda^{2}(2-\rho)=-1\end{equation}
\begin{equation}\label{eq42}-\frac{\delta^{2}}{4}-\frac{\rho}{4}=\frac{\delta}{2\lambda}.\end{equation}
Now, substituting (\ref{eq40}) in (\ref{eq28}) and using (\ref{eq42}) we obtain
\begin{equation}\label{eq43}f=\frac{\chi_{2}+\chi_{3}}{2}-\frac{4-\rho}{2\rho}(\chi_{3}-\chi_{2})\end{equation}
and
\begin{equation}\label{eq44}g=\frac{\chi_{2}+\chi_{3}}{2}-\gamma(\chi_{3}-\chi_{2}),\end{equation}
where $\gamma:=\frac{\delta^{2}}{\rho}+\frac{\delta}{2\lambda}$.
Then multiplying (\ref{eq43}) by $\gamma$ and (\ref{eq44}) by $\frac{4-\rho}{2\rho}$, and adding the identities obtained we get that
$$\frac{4-\rho}{2\rho}g+\gamma f=\frac{4-\rho+2\gamma\rho}{2\rho}\frac{\chi_{2}+\chi_{3}}{2}.$$
Since $f$ and $g$ are linearly independent, we deduce, seeing the identity above and the formulas (\ref{eq43}) and (\ref{eq44}), that $4-\rho+2\gamma\rho\neq0$.
So, that, by putting $\omega:=-\frac{2\gamma\rho}{4-\rho+2\gamma\rho}$ we deduce from the identity above that $$(1+\omega)g-\omega f=\frac{\chi_{2}+\chi_{3}}{2}.$$
Since $\chi_{2}\neq0$, $\chi_{3}\neq0$ and $\chi_{2}\neq\chi_{3}$ we deduce, as in Case A of the proof of \cite [Theorem 3.3]{H.St1}, that the solution occurs in part (4).
\item $\chi_{2}=0$. We go back to (a) by writing $\rho$ instead of $2-\rho$, and $-\lambda$ instead of $\lambda$. So we obtain the solution in part (4).
\item $\chi_{3}=0$. As in (a) we deduce that
\begin{equation*}\xi c\rho-\frac{\rho}{4}-\frac{\eta}{2\lambda}=0\end{equation*}
and
\begin{equation*}2\xi c-\frac{1}{2}-\xi c\rho+\frac{\rho}{4}+\frac{\eta}{2\lambda}=0.\end{equation*}
By adding the two identities above we obtain $\xi c=\frac{1}{4}$ and then $\eta=0$. Then, using (\ref{eq30}) and (\ref{eq29}), we obtain $\delta=-i$ and $\xi=\frac{1}{2}$. So that $c=\frac{1}{2}$ because $\xi c=\frac{1}{4}$.\\
Now, using (\ref{eq28}) we obtain by a small computation
\begin{equation*}f=\frac{\chi_{1}+\chi_{2}}{2}+\frac{\rho-1}{2}\frac{\chi_{1}-\chi_{2}}{2}\end{equation*}
and
\begin{equation*}g=\frac{\chi_{1}+\chi_{2}}{2}+\gamma\frac{\chi_{1}-\chi_{2}}{2},\end{equation*}
where $\gamma=\frac{\rho-1}{2}-\frac{i}{2\lambda}$. Proceeding as (a) we obtain the solution in part (4). This completes the proof of Corollary \ref{cor1}.
\end{enumerate}
\end{enumerate}
\end{proof}

\end{document}